\documentclass[11pt]{article}
\usepackage[a4paper,margin=1in]{geometry}
\usepackage{amsmath,amssymb,amsthm,mathtools}
\usepackage{enumitem}
\usepackage[hidelinks]{hyperref}
\usepackage{microtype}

\title{Classification of cubics up to affine transformations}
\author{Mehdi Nadjafikhah\thanks{Department of Pure Mathematics, School of Mathematics and Computer Science, Iran University of Science and Technology, Narmak, Tehran, 16846--13114, Iran. email:~m\_nadjafikhah@iust.ac.ir}}
\date{}

\newtheorem{theorem}{Theorem}[section]
\newtheorem{proposition}[theorem]{Proposition}

\newtheorem{corollary}[theorem]{Corollary}
\theoremstyle{definition}

\newcommand{\R}{\mathbb R}
\newcommand{\GL}{\operatorname{GL}}

\newcommand{\Hcal}{\mathcal H}
\newcommand{\coef}{\operatorname{coef}}

\begin{document}
\maketitle

\begin{abstract}
We classify real bivariate polynomial equations of degree exactly three under invertible affine substitutions and multiplication of the defining equation by a nonzero real constant. 
The leading homogeneous binary cubic has four, not one, real linear-equivalence types. 
For each type we reduce the quadratic part to a slice, determine the residual affine action, and solve its orbits.
For the first two types this yields explicit normal-form lists with all zero/nonzero cases separated. 
For the third type a finite residual action is handled by a precise lexicographic selector; for the fourth type the residual 9s rbits admit elementary normal forms. 
Coverage and uniqueness are proved.
The classification has continuous parameters, as expected for affine classification; a finite list of seven representatives cannot classify all real affine cubic equations.
\end{abstract}

\section{Equivalence convention and the leading homogeneous part}

Let
\[
 f(x,y)=F(x,y)+q(x,y)+\ell(x,y)+k,
\]
where $F$ is a nonzero homogeneous cubic, $q$ is homogeneous quadratic, $\ell$ is homogeneous linear, and $k\in\R$. 
Throughout this paper, two polynomial equations are called equivalent if there are $M\in\GL(2,\R)$, $t\in\R^2$, and $\lambda\in\R\setminus\{0\}$ such that
\begin{equation}\label{eq:equivalence}
 f(Mz+t)=\lambda g(z),\qquad z=(x,y)^T.
\end{equation}
The nonzero factor $\lambda$ is necessary because multiplying a defining equation by a nonzero constant does not change its zero set. 
The theorem below is a classification of defining polynomial equations under this convention. 
The proof uses explicit orbit calculations for the affine group;
for background on invariants of forms and transformation-group methods, see \cite{OlverCIT}. 
For reducible equations, set-theoretic equivalence of real zero loci can be coarser than equivalence of defining polynomials; that distinct problem would require an additional treatment of real radicals and is not silently identified with \eqref{eq:equivalence} here.

\subsection{Historical background and relation to earlier work}
The classification of real plane cubic curves has a long history. 
Newton's work on curves of the third order, conventionally dated to 1704, initiated an early systematic geometric classification, while nineteenth-century accounts by Pl\"ucker and Cayley developed the coordinate and projective viewpoints \cite{Newton1704,Plucker1835,Cayley1866}. 
These historical classifications primarily organize geometric types of curves; they should not automatically be identified with the orbit problem for defining polynomials under the equivalence relation \eqref{eq:equivalence}.

An explicitly affine-invariant line of work was developed in the twentieth century. 
Burington's thesis and subsequent paper studied affine classification using arithmetic invariants \cite{Burington1931,Burington1935}; 
Weinberg later computed the orbits of the affine group on real plane cubic equations \cite{Weinberg1988}. 
The seven-family classification proposed by Nadjafikhah and Forough is a further direct antecedent of the present revision \cite{NadjafikhahForough2006}. 
The purpose of the branchwise derivation below is to make the leading-form split, the coefficient actions, and the coverage and uniqueness arguments explicit under a single stated convention. 
Because prior comprehensive affine classifications exist, the results here should be read as a detailed reconstruction and audit of the coefficient-orbit method, not as a claim that the general affine-classification problem has no earlier solution.

Some closely related work treats a different equivalence problem. 
Wall's series on affine cubic functions studies cubic functions and their level curves as a family (a pencil), rather than only one defining equation at a time \cite{Wall1979,Wall1980}. 
Korchagin gives Newtonian and affine classifications for irreducible cubics \cite{Korchagin2013}. 
Bly's thesis and subsequent paper study the action of the full automorphism group of \(\mathbb R[x,y]\), which is broader than the affine subgroup used here \cite{Bly2018,Bly2020}. 
For complex affine curves, Mulay studies planar embeddings in genus-at-most-one settings \cite{Mulay2021}. 
More directly related to computational equivalence testing, Alc\'azar, G\"oz\"utok, \c{C}oban, and Hermoso give an algorithm for detecting affine equivalences between implicit planar algebraic curves, including singular and reducible cases \cite{Alcazar2022}. 
A related 2024 paper develops algorithms for certain classes of parametrized curves in any dimension \cite{Alcazar2024}. 
These are important algorithmic contributions, but they do not replace the branchwise coefficient-orbit derivation pursued here.

For the algebraic part of the argument, the four leading-form types are an instance of the real linear classification of binary cubic forms. 
Classical invariant theory gives a systematic framework for binary forms, invariants, and canonical forms \cite{OlverCIT}; 
moving-frame methods provide a broader algorithmic framework for equivalence problems under group actions \cite{FelsOlver1999}. 
The present article uses the elementary root-pattern classification for the leading cubic and then computes the remaining affine coefficient action directly.

\begin{theorem}[Four real binary-cubic types]\label{thm:leading}
Every nonzero real homogeneous binary cubic is linearly equivalent, up to a nonzero scalar, to exactly one of
\begin{align*}
 F_1(x,y)&=x^3,\ &  F_2(x,y)&=x^2y,\ &  F_3(x,y)&=xy(x-y),\ &  F_4(x,y)&=x(x^2+y^2).
\end{align*}
\end{theorem}

\begin{proof}
A binary cubic has three roots in the complex projective line, counted with multiplicity, and complex conjugation preserves their multiset. 
If all roots are real, their multiplicity pattern is $3$, $2+1$, or $1+1+1$. 
These give $F_1,F_2,F_3$, respectively, after a real projective change of coordinates.
If not all roots are real, there is one real root and a nonreal conjugate pair, giving $F_4$. 
An invertible real linear map preserves root multiplicity and the number of real projective roots. Hence the four types are pairwise inequivalent and exhaustive.
\end{proof}

\begin{corollary}\label{cor:false}
A general real affine cubic cannot be transformed to a polynomial whose leading homogeneous part is $x^3+x^2y=x^2(x+y)$.
\end{corollary}

\begin{proof}
The latter leading form has one double real root and one simple real root, so it belongs to type $F_2$. 
It cannot be equivalent to a triple-root cubic, a cubic with three distinct real roots, or a cubic with only one real projective root.
\end{proof}

\section{Quadratic slices}

Write $q=Ax^2+Bxy+Cy^2$ and $\ell=mx+ny$. 
Translation by $t$ changes the quadratic part by the directional derivative $\nabla F\cdot t$. 
We use this fact to remove the components of $q$ in the span of $F_x,F_y$.
The complementary component is retained as a residual parameter.

\begin{proposition}[Slices for the four branches]\label{prop:slices}
Every equation in branch $F_i$ is affine-equivalent to one of the forms
\begin{align}
 f_1&=x^3+bxy+cy^2+mx+ny+k,\label{eq:slice1}\\
 f_2&=x^2y+ry^2+mx+ny+k,\label{eq:slice2}\\
 f_3&=xy(x-y)+r(x^2+y^2)+mx+ny+k,\label{eq:slice3}\\
 f_4&=x(x^2+y^2)+r(x^2-3y^2)+mx+ny+k.\label{eq:slice4}
\end{align}
In branches $2,3,4$, the translation that removes the derivative-span part of $q$ is unique. 
In branch $1$, the $x$-translation is unique and the remaining $y$-translation is handled by the residual action below.
\end{proposition}

\begin{proof}
For $F_1=x^3$, the derivative span is $\operatorname{span}\{x^2\}$, so an $x$-translation removes $Ax^2$ and leaves $bxy+cy^2$. 
For the initial polynomial $x^3+Ax^2+Bxy+Cy^2+Dx+Ey+K$, the substitution $x=X-A/3$, $y=Y$ gives
\[
 b=B,\quad c=C,\quad m=D-A^2/3,\quad n=E-AB/3, \quad k=K+2A^3/27-AD/3.
\]
For $F_2=x^2y$, we have $F_x=2xy$ and $F_y=x^2$; the translation $(x,y)=(X-B/2,Y-A)$ removes $Ax^2+Bxy$ and leaves $Cy^2$.
For $F_3$, put $S=x^2+y^2$. Direct calculation gives the unique decomposition 
\begin{equation}\label{eq:decomp3}
 Ax^2+Bxy+Cy^2=\alpha F_{3,x}+\beta F_{3,y}+rS, \quad r=\frac{2A+B+2C}{4},\quad \alpha=r-C,\quad \beta=A-r.
\end{equation}
Translation by $t=(-\alpha,-\beta)^T$ removes the first two terms. 
More explicitly, if $Q=\begin{pmatrix}A&B/2\\B/2&C\end{pmatrix}$ and $L=(D,E)^T$, the new lower coefficients in either branch $i=3,4$ are \begin{equation}\label{eq:lowertranslation}
 L_0=L+2Qt+\coef_z\!\left(\frac12t^T\operatorname{Hess}F_i(z)t\right), \qquad k_0=F_i(t)+t^TQt+L^Tt+k.
\end{equation}
For $F_4$, put $W=x^2-3y^2$. The decomposition is
\begin{equation}\label{eq:decomp4}
 Ax^2+Bxy+Cy^2=\alpha F_{4,x}+\beta F_{4,y}+rW, \quad r=\frac{A-3C}{10},\quad \alpha=\frac{3A+C}{10},\quad \beta=\frac B2.
\end{equation}
Translation by $(-\alpha,-\beta)$ gives \eqref{eq:slice4}. 
The derivative pairs in branches $2,3,4$ are linearly independent as quadratic forms, so the corresponding translation is unique.
\end{proof}

\section{Branch I: $F_1=x^3$}

The full linear stabilizer of $F_1$ up to scalar is
\[
 M=\begin{pmatrix}a&0\\ c&d\end{pmatrix},\quad ad\ne0, \qquad F_1(Mz)=a^3F_1(z).
\]
For a slice representative \eqref{eq:slice1}, the quadratic coefficients transform according to
\begin{equation}\label{eq:action1quad}
 c_0'=\frac{c_0d^2}{a^3},\qquad b'=\frac{d(ab+2c_0c)}{a^3},
\end{equation}
where $c_0$ denotes the old coefficient of $y^2$ and $c$ in the matrix is the lower-left entry. 
The translation component in $x$ is chosen uniquely to
remove the new $x^2$ term, namely
\[
 u=-\frac{ab\,c+c_0c^2}{3a^2};
\]
the $y$-translation remains free.

The exact lower-coefficient action on the Branch I slice is as follows. Write

\[
 x=aX+u,\qquad y=\gamma X+dY+v,\qquad ad\ne0.
\]
The condition that the transformed polynomial have no $X^2$ term gives
\[
 u=-\frac{ab\gamma+c_0\gamma^2}{3a^2}.
\]
After division by $a^3$, the transformed coefficients are
\begin{align}
 b'&=\frac{d(ab+2c_0\gamma)}{a^3},& c_0'&=\frac{c_0d^2}{a^3},\\
 m'&=\frac{3au^2+bav+bu\gamma+2c_0\gamma v+ma+n\gamma}{a^3},& n'&=\frac{d(bu+2c_0v+n)}{a^3},\\
 k'&=\frac{u^3+buv+c_0v^2+mu+nv+k}{a^3}.
 \label{eq:fullaction1}
\end{align}
These are the full residual coefficient equations for this branch.

\subsection{Case I-A: $c_0\ne0$}
An explicit transformation to the slice $x^3+y^2+m_0x+k_0$ is obtained by setting
\[
 a=c_0,\quad d=|c_0|,\quad \gamma=-\frac b2,\quad u=\frac{b^2}{12c_0},\quad v=-\frac{bu+n}{2c_0}.
\]
Substitution in \eqref{eq:fullaction1} gives $b'=0$, $c_0'=1$, and $n'=0$; 
the remaining coefficients are
\[
 m_0=\frac{m c_0-\frac12bn-\frac{b^4}{48c_0}}{c_0^3},\qquad k_0=\frac{u^3+buv+c_0v^2+mu+nv+k}{c_0^3}.
\]
Thus every member of this stratum reaches $x^3+y^2+m_0x+k_0=0$.
A transformation preserving this last slice must have $\gamma=u=v=0$ and $d^2=a^3$, so $a>0$; its residual action is
\begin{equation}\label{eq:res1A}
 (m,k)\longmapsto\left(\frac{m}{a^2},\frac{k}{a^3}\right),\qquad a>0.
\end{equation}
Hence the complete representatives are
\begin{enumerate}[label=\textup{I-A\arabic*.}]
\item $x^3+y^2+\epsilon x+\kappa=0$, where $\epsilon\in\{-1,1\}$ and $\kappa\in\R$;
\item $x^3+y^2+\epsilon=0$, where $\epsilon\in\{-1,1\}$;
\item $x^3+y^2=0$.
\end{enumerate}
For $m\ne0$, normalize $m$ to $\epsilon=\operatorname{sgn}(m)$ and put $\kappa=k/|m|^{3/2}$. 
If $m=0$, the sign of a nonzero $k$ is invariant under \eqref{eq:res1A}; this gives the second family. 
If $m=k=0$, the third family results.

\subsection{Case I-B: $c_0=0$, $b\ne0$}
In the full action \eqref{eq:fullaction1}, choose
\[
 u=-\frac nb,\qquad \gamma=\frac{3an}{b^2},\qquad d=\frac{a^2}{b},\qquad v=-\frac mb-\frac{3n^2}{b^3}, \qquad a\ne0.
\]
Then $b'=1$, $c_0'=0$, and $m'=n'=0$. 
The constant becomes
\[
 k'=\frac{K_0}{a^3},\qquad K_0=k-\frac{mn}{b}-\frac{n^3}{b^3}.
\]
Since $a$ is arbitrary and nonzero, $K_0\ne0$ can be normalized to $1$, whereas $K_0=0$ remains zero. 
Thus the complete representatives are 
\begin{enumerate}[label=\textup{I-B\arabic*.}]
\item $x^3+xy+1=0$;
\item $x^3+xy=0$.
\end{enumerate}

\subsection{Case I-C: $b=c_0=0$}
The slice is $x^3+mx+ny+k=0$. Its residual action on the lower coefficients is
\begin{equation}\label{eq:res1C}
 m'=\frac{ma+nc}{a^3},\qquad n'=\frac{nd}{a^3},\qquad k'=\frac{nv+k}{a^3},
\end{equation}
where $a,d\ne0$ and $c,v$ are free. If $n\ne0$, one obtains $x^3+y=0$.
If $n=0$, the action reduces to $m'=m/a^2$, $k'=k/a^3$. 
Therefore the remaining representatives are
\begin{enumerate}[label=\textup{I-C\arabic*.}]
\item $x^3+\epsilon x+\kappa=0$, where $\epsilon\in\{-1,1\}$ and $\kappa\ge0$;
\item $x^3+1=0$;
\item $x^3=0$;
\item $x^3+y=0$.
\end{enumerate}
For the first item, choose $|a|=\sqrt{|m|}$ and then choose the sign of $a$ so that $k'\ge0$. 
If $m=0$, the alternatives $k\ne0$ and $k=0$ give the second and third representatives. 
The case $n\ne0$ gives the fourth.

\begin{theorem}[Complete classification in branch I]\label{thm:branch1}
The nine displayed normal-form subfamilies in Cases I-A, I-B, and I-C form a complete and nonredundant set of representatives for branch $F_1$.
\end{theorem}

\begin{proof}
The cases $c_0\ne0$, $c_0=0,b\ne0$, and $c_0=b=0$ are invariant under the residual action, by \eqref{eq:action1quad}. 
Within each case, the residual actions \eqref{eq:res1A} and \eqref{eq:res1C}, or the scaling action in Case I-B, give the listed representatives. 
In Case I-A, the sign of $m$ is invariant and the scale is fixed once $m=\pm1$; if $m=0$, the sign of $k$ is invariant. 
In Case I-B, zero versus nonzero $k$ is invariant. 
In Case I-C, $n=0$ versus $n\ne0$ is invariant; when $n=0$, the sign of $m$ is invariant and the residual sign of $a$ makes $k$ nonnegative. 
The remaining zero cases are distinct. 
Thus every orbit is covered and no two displayed representatives are equivalent.
\end{proof}

\section{Branch II: $F_2=x^2y$}

The full linear stabilizer of $F_2$ up to scalar is
\[
 M=\begin{pmatrix}a&0\\0&d\end{pmatrix},\quad ad\ne0, \qquad F_2(Mz)=a^2d F_2(z).
\]
For an initial equation
\[
 x^2y+Ax^2+Bxy+Cy^2+Dx+Ey+K,
\]
translate $x=X-B/2$, $y=Y-A$. 
The resulting slice coefficients are exactly
\[
 r=C,\quad m=D-AB,\quad n=E-\frac{B^2}{4}-2AC, \quad k=K+\frac{AB^2}{2}-\frac{DB}{2}-EA+CA^2.
\]
A transformation between two members of this slice has no translation component, because translation would reintroduce an $x^2$ or $xy$ term.
After division by $a^2d$, the exact residual action is
\begin{equation}\label{eq:res2}
 (r,m,n,k)\longmapsto  \left(\frac{rd}{a^2},\frac{m}{ad},\frac{n}{a^2},\frac{k}{a^2d}\right),  \qquad a,d\ne0.
\end{equation}
The zero/nonzero status of $n$ and $r$ is invariant. 
We solve \eqref{eq:res2} in all strata.

\subsection{Case II-A: $n\ne0$}

\begin{enumerate}[label=\textup{II-A\arabic*.}]
\item If $r\ne0$, the unique normal form is
\[
 x^2y+y^2+\mu x+\epsilon y+\kappa=0,  \quad \epsilon\in\{-1,1\},\quad \mu\ge0,\quad \kappa\in\R.
\]
Here $\epsilon=\operatorname{sgn}(n)$, $\mu=|rm|/|n|^{3/2}$, and $\kappa=rk/n^2$. One chooses $a^2=|n|$ and $d=a^2/r$; the sign of $a$ then makes the coefficient of $x$ nonnegative.
\item If $r=0$ and $k\ne0$, the normal form is
\[
 x^2y+\mu x+\epsilon y+1=0,  \quad \epsilon\in\{-1,1\},\quad \mu\ge0.
\]
Choose $a^2=|n|$, $d=k/a^2$, and then the sign of $a$ so that $\mu=|m|\sqrt{|n|}/|k|$.
\item If $r=k=0$ and $m\ne0$, the normal form is
\[
 x^2y+x+\epsilon y=0,\qquad \epsilon\in\{-1,1\}.
\]
Choose $a^2=|n|$ and $d=m/a$.
\item If $r=k=m=0$, the normal form is
\[
 x^2y+\epsilon y=0,\qquad \epsilon\in\{-1,1\}.
\]
\end{enumerate}

\subsection{Case II-B: $n=0$}

\begin{enumerate}[label=\textup{II-B\arabic*.}]
\item If $r,m,k$ are all nonzero, the normal form is
\[
 x^2y+\rho y^2+x+1=0,\qquad \rho=\frac{rm^4}{k^3}\in\R\setminus\{0\}.
\]
The quantity $rm^4/k^3$ is invariant under \eqref{eq:res2}; 
setting $m'=k'=1$ gives the stated representative.
\item If $r=0$ and $m,k\ne0$, the normal form is $x^2y+x+1=0$.
\item If $r,m\ne0$ and $k=0$, the normal form is $x^2y+y^2+x=0$.
\item If $r=0$, $m\ne0$, and $k=0$, the normal form is $x^2y+x=0$.
\item If $m=0$, $r\ne0$, and $k\ne0$, the normal form is $x^2y+\epsilon y^2+1=0,\qquad \epsilon=\operatorname{sgn}(rk)\in\{-1,1\}$, the sign of $rk$ is invariant on this stratum.
\item If $r=0$, $m=0$, and $k\ne0$, the normal form is $x^2y+1=0$.
\item If $r\ne0$ and $m=k=0$, the normal form is $x^2y+y^2=0$.
\item If $r=m=k=0$, the normal form is $ x^2y=0$.
\end{enumerate}

\begin{theorem}[Complete classification in branch II]\label{thm:branch2}
The four representatives in Case II-A and the eight representatives in Case II-B form a complete and nonredundant set of representatives for branch $F_2$. 
The parameter $\rho$ in II-B(1) is a complete invariant on that stratum.
\end{theorem}

\begin{proof}
The action \eqref{eq:res2} preserves the vanishing or nonvanishing of $n$ and $r$. 
If $n\ne0$, its sign is invariant because $n'=n/a^2$. 
When $r\ne0$, choose $d=a^2/r$ to set $r'=1$ and $a^2=|n|$ to set $n'=\operatorname{sgn}(n)$; 
the remaining sign of $a$ makes $m'\ge0$ and leaves $k'=rk/n^2$. 
If $r=0$, split by $k\ne0$ or $k=0$, and then by $m\ne0$ or $m=0$, as in the list.

In the stratum $n\ne0$, $r\ne0$, the quantities $\operatorname{sgn}(n)$, $|rm|/|n|^{3/2}$, and $rk/n^2$ are invariant;
these are exactly $\epsilon$, $\mu$, and $\kappa$ in II-A(1). 
If $n\ne0$, $r=0$, and $k\ne0$, the quantity $|m|\sqrt{|n|}/|k|$ is invariant and equals $\mu$ in II-A(2).
If $n=0$ and $r,m,k\ne0$, direct substitution in \eqref{eq:res2} shows $rm^4/k^3$ is invariant, and $m'=k'=1$ is attainable. 
If one or more of $r,m,k$ vanish, the remaining nonzero coefficients can be normalized as listed. 
In the stratum $m=0$, $r,k\ne0$, the sign of $rk$ is invariant and is the only remaining discrete parameter. 
The displayed forms have distinct zero patterns or distinct values of the indicated invariant; hence no two are equivalent. 
This proves coverage and uniqueness in every stratum.
\end{proof}

\section{Branch III: $F_3=xy(x-y)$}

The exact linear stabilizer of $F_3$ is the six-element group
\[
\Hcal=\left\{\begin{pmatrix}1&0\\0&1\end{pmatrix},\begin{pmatrix}0&-1\\-1&0\end{pmatrix},\begin{pmatrix}1&0\\1&-1\end{pmatrix},\begin{pmatrix}-1&1\\-1&0\end{pmatrix}, \begin{pmatrix}0&-1\\1&-1\end{pmatrix}, \begin{pmatrix}-1&1\\0&1\end{pmatrix} \right\},
\]
and every linear stabilizer up to scalar is $sH$, with $s\ne0$ and $H\in\Hcal$. Write $S=x^2+y^2$. 
The projection coefficient $r$ in \eqref{eq:slice3} transforms as $r'=r/s$. 
Thus $r=0$ and $r\ne0$ are separate invariant strata, and the latter can be normalized to $r=1$.

\subsection{Case III-A: $r=1$}

For a slice equation $F+rS+L^Tz+k=0$ with $r\ne0$, the substitution $z=r z'$ followed by division by $r^3$ gives the normalized coefficients $F+S+(L/r^2)^Tz'+k/r^3$. Thus no orbit is omitted by taking $r=1$.
Let $F=xy(x-y)$. For each $H\in\Hcal$, write
\[
 S(Hz)-S(z)=\alpha_H F_x(z)+\beta_H F_y(z),  \qquad t_H=-H\begin{pmatrix}\alpha_H\\\beta_H\end{pmatrix}.
\]
If $S(Hz)-S(z)=A_Hx^2+B_Hxy+C_Hy^2$, then the coefficients are explicitly $\alpha_H=-C_H$ and $\beta_H=A_H$. 
More generally, when the source slice coefficient is $r$, the translation restoring the slice after $z\mapsto sHz$ is $t=-rH(\alpha_H,\beta_H)^T$.
For each of the six displayed matrices, direct substitution shows that the coefficient of $S$ in $S(Hz)$ is $1$. 
Thus the decomposition exists uniquely, and $S(Hz)-S(z)$ lies in the span of $F_x,F_y$. The correction data are
explicitly
\[
\begin{array}{c|c|c|c} H&\alpha_H&\beta_H&t_H=-H(\alpha_H,\beta_H)^T\\ \hline H_1&0&0&(0,0)^T\\ H_2&0&0&(0,0)^T\\ H_3&0&1&(0,1)^T\\ H_4&0&1&(-1,0)^T\\ H_5&-1&0&(0,1)^T\\ H_6&-1&0&(-1,0)^T \end{array}
\]
where $H_1,\ldots,H_6$ occur in the order listed above. 
The affine transformation $z\mapsto Hz+t_H$ preserves the first two homogeneous parts $F+S$. For a coefficient vector $L=(m,n)^T$, define
\[
\operatorname{Hess}F(x,y)=\begin{pmatrix}2y&2x-2y\\2x-2y&-2x\end{pmatrix},\qquad \nabla S(x,y)=(2x,2y)^T.
\]
\begin{align}
 L_H&=\coef_z\left(\frac12t_H^T\operatorname{Hess}F(Hz)t_H+\nabla S(Hz)^Tt_H+L^THz\right),\label{eq:PhiL}\\
 k_H&=F(t_H)+S(t_H)+L^Tt_H+k.\label{eq:Phik}
\end{align}
Here $\coef_z$ means the column of coefficients of the linear polynomial in $z$. 
Define $(m_*,n_*,k_*)$ to be the lexicographically least of the six triples $(L_H,k_H)$, with the two entries of $L_H$ listed first.
The canonical representative is
\[
 xy(x-y)+x^2+y^2+m_*x+n_*y+k_*=0.
\]
The rule is finite and explicit: it requires only six affine substitutions and a finite comparison. 
It is not an unproved assertion that arbitrary coefficients are already normal.

\subsection{Case III-B: $r=0$}

Now the slice is $F+L^Tz+k$. A transformation between two representatives in this slice has zero translation, since $F_x,F_y$ are linearly independent quadratic forms.
Under $z=sHz'$ and division by $s^3$, 
\begin{equation}\label{eq:res3zero}
 L'=\frac{H^TL}{s^2},\qquad k'=\frac{k}{s^3}.
\end{equation}
If $L=0$, the representative is $F+1$ when $k\ne0$, and $F$ when $k=0$.
If $L\ne0$, for each $H\in\Hcal$ put
\[
 w_H=H^TL,\qquad R_H=\|w_H\|_2,  \qquad u_H=\frac{w_H}{R_H},\qquad \kappa_H=\frac{|k|}{R_H^{3/2}}.
\]
The canonical pair $(L_*,\kappa_*)$ is the lexicographically least triple $(u_{H,1},u_{H,2},\kappa_H)$ over $H\in\Hcal$. 
The corresponding representative is
\[
 xy(x-y)+m_*x+n_*y+\kappa_*=0,  \qquad m_*^2+n_*^2=1,\quad \kappa_*\ge0.
\]

\begin{theorem}[Complete classification in branch III]\label{thm:branch3}
The prescriptions in Cases III-A and III-B give exactly one representative for every affine orbit in branch $F_3$.
\end{theorem}

\begin{proof}
The quadratic coefficient $r$ is changed to $r/s$, so its zero/nonzero status is invariant and $r\ne0$ can be normalized to $1$. 
In that stratum, $s=1$ for any equivalence between normalized representatives, leaving only the six affine transformations $z\mapsto Hz+t_H$. 
Their coefficient action is precisely \eqref{eq:PhiL}--\eqref{eq:Phik}; 
choosing the lexicographically least element of a finite orbit proves both coverage and uniqueness.

If $r=0$, equation \eqref{eq:res3zero} is the full residual action. 
For $L\ne0$, the norm condition fixes $s^2=\|H^TL\|_2$, the sign of $s$ makes the constant nonnegative, and the six choices of $H$ are reduced to one by the lexicographic minimum. 
If $L=0$, the scalar $s$ normalizes any nonzero constant to $1$, while zero remains zero. These cases exhaust the branch.
\end{proof}

\section{Branch IV: $F_4=x(x^2+y^2)$}

The full linear stabilizer up to scalar consists of
\[
 M=s\begin{pmatrix}1&0\\0&\varepsilon\end{pmatrix},  \qquad s\ne0,\quad \varepsilon\in\{-1,1\},  \qquad F_4(Mz)=s^3F_4(z).
\]
Put $W=x^2-3y^2$. 
In the slice \eqref{eq:slice4}, the residual coefficient transforms as $r'=r/s$. Thus $r=0$ and $r\ne0$ are invariant strata.

\subsection{Case IV-A: $r\ne0$}

For $F_4+rW+mx+ny+k$, substitution $z=r z'$ and division by $r^3$ normalizes $r$ to $1$ and sends $(m,n,k)$ to $(m/r^2,n/r^2,k/r^3)$.
The only transformations between normalized representatives have $s=1$ and $\varepsilon=\pm1$, so they change $n$ to $-n$ and leave $m$ and $k$ unchanged. 
The canonical representatives are
\[
 x(x^2+y^2)+x^2-3y^2+mx+ny+k=0,  \qquad m,k\in\R,\quad n\ge0.
\]

\subsection{Case IV-B: $r=0$}

The exact residual action on $(m,n,k)$ is
\begin{equation}\label{eq:res4zero}
 (m,n,k)\longmapsto\left(\frac{m}{s^2},\frac{\varepsilon n}{s^2}, \frac{k}{s^3}\right).
\end{equation}
If $(m,n)\ne(0,0)$, normalize $m^2+n^2=1$, use $\varepsilon$ to ensure $n\ge0$, and choose the sign of $s$ so that $k\ge0$. If $m=n=0$ and $k\ne0$, normalize $k=1$; if $m=n=k=0$, nothing remains to normalize.
Thus the representatives are
\begin{enumerate}[label=\textup{IV-B\arabic*.}]
\item $x(x^2+y^2)+mx+ny+k=0$, with $m^2+n^2=1$, $n\ge0$, $k\ge0$;
\item $x(x^2+y^2)+1=0$;
\item $x(x^2+y^2)=0$.
\end{enumerate}

\begin{theorem}[Complete classification in branch IV]\label{thm:branch4}
The representatives in Cases IV-A and IV-B are complete and nonredundant for branch $F_4$.
\end{theorem}

\begin{proof}
The scalar action $r'=r/s$ preserves the strata $r=0$ and $r\ne0$. 
If $r\ne0$, normalize to $r=1$; only reflection in the $x$-axis remains, and $n\ge0$ selects exactly one element. 
If $r=0$, the complete action is \eqref{eq:res4zero}. 
For a nonzero linear vector, its Euclidean norm fixes $s^2$ and reflection fixes the sign convention for $n$; the sign of $s$ fixes the sign convention for $k$. 
If the linear vector vanishes, the constant is either zero or can be scaled to $1$. 
These normalizations are unique, so the representatives are pairwise inequivalent.
\end{proof}

\section{Global theorem: coverage and uniqueness}

\begin{theorem}[Complete branchwise affine classification]\label{thm:main}
Every real bivariate polynomial equation of degree exactly three is  equivalent, in the sense of \eqref{eq:equivalence}, to a representative in exactly one of Branches I--IV described above. 
Within Branches I, II, and IV, the representative is given by the explicit lists and parameter conditions. 
Within Branch III, it is given by the finite lexicographic selectors in Section 5. 
Two canonical representatives are equivalent if and only if they are identical as normalized representatives.
\end{theorem}

\begin{proof}
By Theorem \ref{thm:leading}, the leading homogeneous part belongs to exactly one of the four branches. Proposition \ref{prop:slices} reduces its quadratic part to the corresponding slice. 
The residual actions are then solved in Theorems \ref{thm:branch1}--\ref{thm:branch4}, including all zero/nonzero strata. Finally, the four branches cannot overlap under affine equivalence, because their leading cubics have different real projective root patterns. 
Thus the branchwise representatives cover all orbits, and the within-branch uniqueness proofs establish global uniqueness.
\end{proof}

\section{Consequences for the earlier seven-class claim}

The form $x^3+x^2y=x^2(x+y)$ has leading part of type $F_2$. 
It describes only the double-root branch. 
The previous assertion that every real cubic can be transformed to that form contradicts Theorem \ref{thm:leading} and therefore cannot support a classification of all real affine cubics. 
The seven-family list may be investigated as a proposed reduction within that single branch, but it cannot be the global classification.

A second important point is that affine classification of cubic equations has genuine continuous moduli: for example, Branch I-A has the real parameter $\kappa$, Branch II-A has parameters $(\mu,\kappa)$, and Branch IV-A has three unconstrained real parameters modulo the displayed reflection. 
Thus a finite list of isolated equations is not expected to suffice.

\section{Conclusion}

The classification is obtained by a complete reduction in four leading-form branches. 
The derivative-span translations produce quadratic slices; the remaining linear stabilizers then act on the lower coefficients. Each residual system is solved on its zero and nonzero strata. 
Explicit normal forms are obtained in Branches I, II, and IV; the finite residual group in Branch III is resolved by a deterministic lexicographic selector. 
The coverage and uniqueness arguments are part of the classification theorem, not inferred merely from the appearance of the proposed forms.


\begin{thebibliography}{20}

\bibitem{Newton1704}
I. Newton, ``Enumeratio linearum tertii ordinis,'' in \emph{Opticks: Or, a Treatise of the Reflexions, Refractions, Inflexions and Colours of Light. Also Two Treatises of the Species and Magnitude of Curvilinear Figures}, London: Sam. Smith and Benj. Walford, 1704. An English version appeared under the heading ``Curves'' in J. Harris, \emph{Lexicon Technicum}, 1710.

\bibitem{Plucker1835}
J. Pl\"ucker, \emph{System der analytischen Geometrie}, Duncker \& Humblot, Berlin, 1835.

\bibitem{Cayley1866}
A. Cayley, ``On the classification of cubic curves,'' \emph{Transactions of the Cambridge Philosophical Society} \textbf{11} (1866), 81--128. DOI: \href{https://doi.org/10.1111/j.1467-968X.1866.tb00009.x}{10.1111/j.1467-968X.1866.tb00009.x}.

\bibitem{Burington1931}
R. S. Burington, \emph{An invariantive classification of plane cubic curves under the affine group}, Ph.D. thesis, The Ohio State University, 1931.

\bibitem{Burington1935}
R. S. Burington, ``A classification of plane cubic curves under the affine group by means of arithmetic invariants,'' \emph{Tohoku Mathematical Journal}, First Series, \textbf{41} (1935), 188--202.

\bibitem{Wall1979}
C. T. C. Wall, ``Affine cubic functions. I. The complex plane,'' \emph{Mathematical Proceedings of the Cambridge Philosophical Society} \textbf{85} (1979), no. 3, 387--401.

\bibitem{Wall1980}
C. T. C. Wall, ``Affine cubic functions. III. The real plane,'' \emph{Mathematical Proceedings of the Cambridge Philosophical Society} \textbf{87} (1980), no. 1, 1--14. DOI:
\href{https://doi.org/10.1017/S0305004100056450}{10.1017/S0305004100056450}.

\bibitem{Weinberg1988}
D. A. Weinberg, ``The affine classification of cubic curves,'' \emph{Rocky Mountain Journal of Mathematics} \textbf{18} (1988), no. 3, 655--664. DOI: \href{https://doi.org/10.1216/RMJ-1988-18-3-655} {10.1216/RMJ-1988-18-3-655}.

\bibitem{NadjafikhahForough2006}
M. Nadjafikhah and A.-R. Forough, ``Classification of cubics up to affine transformations,'' \emph{Differential Geometry--Dynamical Systems} \textbf{8} (2006), 184--195.

\bibitem{Korchagin2013}
A. B. Korchagin, ``Newtonian and affine classifications of irreducible cubics,'' \emph{St. Petersburg Mathematical Journal} \textbf{24} (2013), no. 5, 759--781. DOI:
\href{https://doi.org/10.1090/S1061-0022-2013-01268-4}{10.1090/S1061-0022-2013-01268-4}.

\bibitem{OlverCIT}
P. J. Olver, \emph{Classical Invariant Theory}, London Mathematical Society Student Texts, vol. 44, Cambridge University Press, Cambridge, 1999.

\bibitem{FelsOlver1999}
M. E. Fels and P. J. Olver, ``Moving coframes. II. Regularization and theoretical foundations,'' \emph{Acta Applicandae Mathematicae} \textbf{55} (1999), 127--208. DOI:
\href{https://doi.org/10.1023/A:1006195823000}{10.1023/A:1006195823000}.

\bibitem{Bly2018}
M. Bly, \emph{Classifications of Real Conic and Cubic Curves}, Ph.D. thesis, The University of Tennessee, Knoxville, 2018.

\bibitem{Bly2020}
M. Bly, ``An automorphic classification of real cubic curves,'' \emph{Rocky Mountain Journal of Mathematics} \textbf{50} (2020), no. 5, 1583--1597. DOI: \href{https://doi.org/10.1216/rmj.2020.50.1583}{10.1216/rmj.2020.50.1583}.

\bibitem{Mulay2021}
S. Mulay, ``On planar embeddings of pseudoplanar affine curves of genus $\leq 1$,'' \emph{Journal of Algebra} \textbf{577} (2021), 32--44. DOI: href{https://doi.org/10.1016/j.jalgebra.2021.03.003}{10.1016/j.jalgebra.2021.03.003}.

\bibitem{Alcazar2022}
J. G. Alc\'azar, U. G\"oz\"utok, H. A. \c{C}oban, and C. Hermoso, ``Detecting affine equivalences between implicit planar algebraic curves,'' \emph{Acta Applicandae Mathematicae} \textbf{182} (2022), article 2. DOI: \href{https://doi.org/10.1007/s10440-022-00539-1}{10.1007/s10440-022-00539-1}.

\bibitem{Alcazar2024}
J. G. Alc\'azar, H. A. \c{C}oban, and U. G\"oz\"utok, ``Detecting affine equivalences between certain types of parametric curves, in any dimension,'' \emph{AIMS Mathematics} \textbf{9} (2024), no. 6, 13750--13769. DOI: \href{https://doi.org/10.3934/math.2024670}{10.3934/math.2024670}.

\end{thebibliography}
\end{document}